\documentclass[final,1p,times]{elsarticle}

\usepackage[colorlinks=true]{hyperref}
\usepackage{amsfonts}
\usepackage[all]{xy}
\usepackage[mathscr]{eucal}
\usepackage{extarrows}
\usepackage{amsthm,amsmath,amssymb,mathrsfs,amsfonts,graphicx,graphics,latexsym,exscale,cmmib57,dsfont,bbm,amscd}
\usepackage{color,soul}
\newtheorem{thm}{Theorem}[section]

\newtheorem{lem}[thm]{Lemma}

\newtheorem{prop}[thm]{Proposition}
\theoremstyle{definition}
\newtheorem{remark}[thm]{Remark}
\theoremstyle{definition}
\newtheorem{defn}[thm]{Definition}

\numberwithin{equation}{section}
\journal{xxx}

\begin{document}

\begin{frontmatter}
\title{Relative entropy dimensions for amenable group actions}

\author{Zubiao Xiao\corref{cor1}}
\ead{xzb2020@fzu.edu.cn}
\address{School of Mathematics and Statistics, Fuzhou University, Fuzhou 350116, People's Republic of China}
\author{Zhengyu Yin}
\ead{mg1921038@smail.nju.edu.cn}
\address{Department of Mathematics, Nanjing University, Nanjing 210093, People's Republic of China}

\begin{abstract}
We study the topological complexities of relative entropy zero extensions acted by countable-infinite amenable groups. Firstly, for a given F{\o}lner sequence $\{F_n\}_{n=0}^{+\infty}$, we define respectively the relative entropy dimensions and the dimensions of the relative  entropy generating sets to characterize the sub-exponential growth of the relative topological complexity. Meanwhile, we investigate the relations among them. Secondly, we introduce the notion of a relative dimension set. Moreover, using it, we discuss the disjointness between the relative entropy zero extensions which generalizes the results of Dou, Huang and Park[Trans. Amer. Math. Soc. 363(2) (2011), 659-680].
\end{abstract}

\begin{keyword}
amenable groups $\cdot$ relative entropy dimension $\cdot$ relative dimension sets

\medskip
\MSC[2020]  37B99 $\cdot$ 54H15
\end{keyword}
\end{frontmatter}

\section{Introduction}\label{sec1}
A dynamical system for a group action is usually written by a pair $(X,G)$ where $X$ is a set (called a phase space) and $G$ is a group which acts continuously on $X$ (called the acted group). In the study of dynamical system, entropy is an important tool to characterize the dynamical behavior. Kolomogorov and Sinai developed the measure-theoretic entropy of $\mathbb{Z}$-actions based on Shannon's information theory in 1959. Topological entropy was first introduced by Adler, Konheim and McAndrew and defined by Bowen later on a metric space in 1973. Topological entropy measures the maximal exponential growth rate of orbits for an arbitrary topological dynamical system.

Although systems with positive entropy are much more complicated than those with zero entropy, zero entropy systems own various levels of complexity, and recently have been discussed in \cite{Bu,Ca,DHP,DHP1,FP,HPY,KCML,Zhou}. Those authors adopted various methods to classify zero entropy dynamical systems. Carvalho \cite{Ca} introduced the notion of entropy dimension to distinguish the zero entropy systems and obtained some basic properties of entropy dimension. Ferenczi and Park \cite{FP} proposed the entropy dimension for the action of $\mathbb{Z}$ or $\mathbb{Z}^d$ on a probability space. Dou, Huang and Park in \cite{DHP} introduced the notion of the dimension of subset of $\mathbb{Z}$ with density $0$. And then they used the dimension of a special class of sequences which were called entropy generating sequences to measure the complexity of a system and showed that the topological entropy dimension can be computed through the dimensions of entropy generating sequences. Zhou \cite{Zhou} defined the corresponding notions for $\mathbb{Z}$-actions and studied their properties. In this paper, we would like to redefine the relevant notions in relative setting for group actions through introducing a new notion which we call relative "relative entropy generating sets". Using this concept, we can define the relative dimensions of entropy generating sets and get the relationships among different relative dimensions. We can show that the relative upper entropy dimension of the extension is the supremum of the dimensions of the relative upper entropy generating sets (see Theorem \ref{t1}).

Inspired by the theory of prime numbers, Furstenberg \cite{Fu} first introduced the concept of disjointness to characterize the difference of dynamics between two systems. Two important results about disjointness are that in measurable dynamics $K$-mixing systems are disjoint from ergodic zero entropy systems and that weak mixing systems are disjoint from group rotations. In the case of topological settings, these properties are explored in \cite{Bl,HPY,HY}. Dou, Huang and Park in \cite{DHP} introduced the notion of the dimension set of a zero entropy topological system to measure the various levels of topological complexity of subexponential growth rate. They investigated the property of disjointness in zero entropy systems through the dimension set and proved that under the condition of one system's minimality, two systems with disjoint dimension sets are disjoint. This is a refinement and also a generalization of the result in \cite{Bl} that uniformly positive entropy systems are disjoint from minimal entropy zero systems. Based on the above results, Zhou \cite{Zhou} introduced the notion of relative dimension tuples and the relative dimension set and proved that two extensions with disjoint relative dimension sets for all orders are disjoint over the same system under some conditions. This can be also regarded as a generalization of the results in \cite{HYZ1} that an open extension with relative uniformly positive entropy of all orders is disjoint from any minimal extension with relative zero entropy. We would like to redefine the notion of relative dimension tuples and the relative dimension set for group actions and show that an open extension is disjoint from a minimal extension if they have disjoint $n$-th relative dimension sets, for any $n\geq 2$ (see Theorem \ref{t2}).

The paper is organized as follows. In Section \ref{sec2}, we give the definition of amenable groups and some basic concepts of dynamical systems. In Section \ref{sec3}, we define the relative entropy dimension of an extension for an amenable group action and investigate the relevant properties. In Section \ref{sec4}, we consider the dimensions of the relative entropy generating sets. In Section \ref{sec5}, we study the interrelations among the previous defined dimensions. In Section \ref{sec6}, we give the notions of relative dimension tuples, dimension sets and study inheritance and lifting properties of the uniform relative entropy dimension extension. In Section \ref{sec7}, we prove the disjointness theorem between the extensions with disjoint relative dimension sets.
\section{Preliminaries}\label{sec2}
\subsection{Amenable groups}\label{sec2.1}
Let $G$ be a countable discrete infinite group. Denote by $\mathcal{F}(G)$ the set of all finite non-empty subsets of $G$. For $K,F\in \mathcal{F}(G)$, we write
\[KF\bigtriangleup F=(KF\setminus F)\cup (F\setminus KF),\]
where $KF=\{kf:k\in K,f\in F\}$. The group $G$ is called \textit{amenable} if for any $K\in \mathcal{F}(G)$ and $\delta>0$, there exists $F\in \mathcal{F}(G)$ such that \[\frac{|KF\bigtriangleup F|}{|F|}<\delta,\] where $|\cdot|$ is the counting measure on $G$. Such a set $F$ is called $(K,\delta)$-\textit{invariant}. A sequence $\{F_n\}_{n\in \mathbb{N}}\subseteq \mathcal{F}(G)$ is called \textit{a F${\o}$lner sequence} if for every $K\in \mathcal{F}(G)$ and $\delta>0$, for all large enough $n$ we have that $F_n$ is $(K,\delta)$-invariant . A amenable group $G$ is amenable if and only if $G$ admits a F${\o}$lner sequence $\{F_n\}_{n\in \mathbb{N}}$. The class of amenable groups contains in particular all finite groups, all abelian groups and more generally, all solvable groups. It is closed under the operations of taking subgroups, taking quotients, taking extensions and taking inductive limits. For more details and properties of the amenable group, one can refer to \cite[Chaper 4]{CC}.
\subsection{$G$-systems and related concepts}\label{sec2.2}
A $G$-system is a pair $(X,G)$ where $X$ is a compact metric space and $G$ is a group which acts continuously on $X$. Suppose that $e$ is the identity element of $G$. Each element $g\in G$ will be regarded as a homeomorphic action from $X$ to itself if there is no confusion. When $X$ is a set consisting of a singer point, we call $(X,G)$ a trivial system. For a system $(X,G)$ and a positive integer $k\geq 2$, \textit{a $k$-product $G$-system of $(X,G)$} is denoted by $(X^k, G)$ and we assume that $g(x_1,\cdots,x_k)=(gx_1,\cdots,gx_k)$, for every $(x_1,\cdots,x_k)\in X^k$ and $g\in G$.

Let $(X,G)$ be a $G$-system and $x\in X$. For a subset $F\subseteq G$, we denote the $F$-orbit of $x$ by $Fx=\{gx:g\in F\}$. We will call the orbit of $x$ instead of $G$-orbit of $x$ if there is no confusion. The subset $K$ of $X$ is said to be \textit{$G$-invariant} if $GK=K$ (equivalently $GK\subseteq K$). Thus a set is invariant if and only if it is a union of orbits.
\section{Relative entropy dimensions}\label{sec3}
Let $X$ be a compact metric space and $G$ a countable discrete infinite amenable group. Let $(X,G)$ be a $G$-system. Without loss of generality, we assume that $G$ admits a strictly increasing F${\o}$lner sequence $\{F_n\}_{n=0}^{+\infty}$, where $e\in F_0$. If not, we can take a subsequence denoted by $\{F'_n\}_{n=0}^{+\infty}$, Then we choose some sequence $\{n_j\}_{j=1}^\infty$ with $n_1<n_2<\cdots$ to construct $\{F^*_n\}$ as follows: $F_1^*=\{e\}\cup F'_{n_1}$, $F_2^*=F_1^*\cup F'_{n_2}$, $F_3^*=F_2^*\cup F'_{n_3}$, $\cdots$, $F_k^*=F_{k-1}^*\cup F'_{n_k}$, $\cdots$. It is easy to check that $\{e\}\subseteq F_1^*\subseteq F_2^*\subseteq \cdots$ and $\{F_n^*\}$ forms a F${\o}$lner sequence of $G$.

Before we introduce the concepts of relative entropy dimensions of an extension, we recall some definitions. Given a $G$-system $(X,G)$, denote by $\mathcal{C}_X$ the set of finite covers of $X$ and $\mathcal{C}^o_X$ the set of finite open covers of $X$. Given two covers $\mathcal{U},\mathcal{V}\in \mathcal{C}_X$, we say that $\mathcal{U}$ is finer than $\mathcal{V}$ (denoted by $\mathcal{U}\succeq\mathcal{V}$) if for every $U\in \mathcal{U}$, there exists $V\in \mathcal{V}$ such that $U\subseteq V$. Let $\mathcal{U}\vee\mathcal{V}=\{U\cap V:U\in \mathcal{U},V\in \mathcal{V}\}$. Clearly, $\mathcal{U}\vee\mathcal{V}\succeq \mathcal{U}$ and $\mathcal{U}\vee\mathcal{V}\succeq \mathcal{V}$.

Let $(X,G)$ and $(Y,G)$ be two $G$-systems. We call $(X,G)$ is an extension of $(Y,G)$ if there exists a continuous surjective map $\pi:(X,G)\rightarrow (Y,G)$ such that $\pi\circ g=g\circ \pi$ for all $g\in G$. The map $\pi$ is called a factor map from $X$ to $Y$.

Let $\pi:(X,G)\rightarrow (Y,G)$ be a factor map and $\mathcal{U}\in \mathcal{C}_X^o$. For a closed set $E\subseteq X$, we denote
\[
\mathcal{N}(\mathcal{U},E)=\min|\{U_i\in\mathcal{U}:\bigcup\limits_iU_i\supseteq E\}|
\]
and $\mathcal{N}(\mathcal{U}|\pi)=\sup\limits_{y\in Y}\mathcal{N}(\mathcal{U},\pi^{-1}(y))$. For $\alpha\geq 0$, we define
\[
\overline{h}(G,\mathcal{U},\alpha|\pi)=\limsup\limits_{n\rightarrow \infty}\frac{\log \mathcal{N}(\bigvee_{g\in F_n}g^{-1}\mathcal{U}|\pi)}{|F_n|^\alpha},
\]
\[
\underline{h}(G,\mathcal{U},\alpha|\pi)=\liminf\limits_{n\rightarrow \infty}\frac{\log \mathcal{N}(\bigvee_{g\in F_n}g^{-1}\mathcal{U}|\pi)}{|F_n|^\alpha}.
\]
It is clear that $\overline{h}(G,\mathcal{U},\alpha|\pi)$ does not decrease as $\alpha$ decreases, and $\overline{h}(G,\mathcal{U},\alpha|\pi)\notin\{0,+\infty\}$ for at most one $\alpha\geq 0$. We define the relative upper entropy dimension of $\mathcal{U}$ relevant to $\pi$ by
\[
\overline{D}(G,\mathcal{U}|\pi)=\inf\{\alpha\geq0:\overline{h}(G,\mathcal{U},\alpha|\pi)=0\}=\sup\{\alpha\geq0:\overline{h}(G,\mathcal{U},\alpha|\pi)=+\infty\}.
\]
Similarly, $\underline{h}(G,\mathcal{U},\alpha|\pi)$ does not decrease as $\alpha$ decreases, and $\underline{h}(G,\mathcal{U},\alpha|\pi)\notin\{0,+\infty\}$ for at most one $\alpha\geq 0$. We define the relative lower entropy dimension of $\mathcal{U}$ relevant to $\pi$ by
\[
\underline{D}(G,\mathcal{U}|\pi)=\inf\{\alpha\geq0:\underline{h}(G,\mathcal{U},\alpha|\pi)=0\}=\sup\{\alpha\geq0:\underline{h}(G,\mathcal{U},\alpha|\pi)=+\infty\}.
\]
If $\overline{D}(G,\mathcal{U}|\pi)=\underline{D}(G,\mathcal{U}|\pi)=\alpha$, then we say $\mathcal{U}$ has relative entropy dimension $\alpha$. Clearly $0\leq \underline{D}(G,\mathcal{U}|\pi)\leq \overline{D}(G,\mathcal{U}|\pi)\leq 1$ and if \[h(G,\mathcal{U},1|\pi)=\lim\limits_{n\rightarrow \infty}\frac{\log \mathcal{N}(\bigvee_{g\in F_n}g^{-1}\mathcal{U}|\pi)}{|F_n|}>0,\]
i.e, $\mathcal{U}$ has positive relative topological entropy, then the relative entropy dimension of $\mathcal{U}$ is equal to 1.
\begin{defn}\label{d1}
Let $\pi:(X,G)\rightarrow (Y,G)$ be a factor map between $G$-systems. The relative upper (resp. lower) entropy dimension of $(X,G)$ relevant to $\pi$ is
\[
\overline{D}(X,G|\pi)=\sup\{\overline{D}(G,\mathcal{U}|\pi):\mathcal{U}\in\mathcal{C}_X^o\}\quad
(\mbox{resp. }\underline{D}(X,G|\pi)=\sup\{\underline{D}(G,\mathcal{U}|\pi):\mathcal{U}\in\mathcal{C}_X^o\})
\]
\end{defn}
It is clear that $0\leq \underline{D}(X,G|\pi)\leq \overline{D}(X,G|\pi)\leq 1$. When $\overline{D}(X,G|\pi)=\underline{D}(X,G|\pi)$, we call the common value the relative entropy dimension of $(X,G)$ relevant to $\pi$, denoted by $D(X,G|\pi)$. When $(Y,G)$ is a trivial system, "relative" can be removed and the dimension is the classical entropy dimension (see \cite{DHP}), and in this case we shall omit the restriction on $\pi$. Relevant results for measure-theoretic settings can be found in recent work by D. Dou, W. Huang and K. K. Park in \cite{DHP1}.
\begin{remark}

It is known that the definition of entropy for a amenable group action is independent on the choice of F${\o}$lner sequences. However, we consider the definition of entropy dimension of zero entropy system for a amenable group action and this is not the case. It is heavily dependent on the choice of F${\o}$lner sequences. We give a simple example as follows: In \cite[Example 2.8]{DHP}, a shift system $(Y_1,\sigma_1)$ is obtained with the upper entropy dimension $\overline{D}(Y_1,\sigma_1)=\overline{D}(\sigma_1,\mathcal{U}_1)=\frac{1}{2}$ and the lower entropy dimension $\underline{D}(Y_1,\sigma_1)=\underline{D}(\sigma_1,\mathcal{U}_1)=0$. Let $0<\alpha<\frac{1}{2}$. By the definitions of the upper and lower entropy dimension, we have that
\[
\limsup\limits_{n\rightarrow \infty}\frac{\log \mathcal{N}(\bigvee_{i=0}^{n-1}\sigma_1^{-i}\mathcal{U}_1)}{n^\alpha}=+\infty\mbox{ and }\liminf\limits_{n\rightarrow \infty}\frac{\log \mathcal{N}(\bigvee_{i=0}^{n-1}\sigma_1^{-i}\mathcal{U}_1)}{n^\alpha}=0.
\]
By the definition of $\liminf$, we can take a subsequence $\{n_k\}_{k=1}^{\infty}$ of $\{n\}$ such that \[\lim\limits_{k\rightarrow \infty}\frac{\log \mathcal{N}(\bigvee_{i=0}^{n_k-1}\sigma_1^{-i}\mathcal{U}_1)}{n_k^\alpha}=0.\]
We set $F'_k=[0,n_k-1]$. Then $\{F'_k\}_{k=1}^\infty$ is a F${\o}$lner sequence. But in this case, \[\overline{D}(\sigma_1,\mathcal{U}_1,\{[0,n_k-1]\}_{k=1}^\infty)\leq\alpha<\frac{1}{2}.\]
Therefore, $\overline{D}(\sigma_1,\mathcal{U}_1,\{[0,n_k-1]\}_{k=1}^\infty)\neq \overline{D}(\sigma_1,\mathcal{U}_1,\{[0,n-1]\}_{n=1}^\infty)$.
\end{remark}
The following two propositions are basic properties of relative entropy dimensions.
\begin{prop}\label{p1}
Let $\pi:(X,G)\rightarrow (Y,G)$ be a factor map between $G$-systems and $\mathcal{U},\mathcal{V}\in \mathcal{C}_X^o$.
\begin{enumerate}[(1)]
  \item If $\mathcal{U}\preceq\mathcal{V}$, then $\overline{D}(G,\mathcal{U}|\pi)\leq\overline{D}(G,\mathcal{V}|\pi)$ and $\underline{D}(G,\mathcal{U}|\pi)\leq\underline{D}(G,\mathcal{V}|\pi)$.
  \item If $G$ is an abelian group, then $\overline{D}(G,\mathcal{U}|\pi)=\overline{D}(G,g\mathcal{U}|\pi)$ and $\underline{D}(G,\mathcal{U}|\pi)=\underline{D}(G,g\mathcal{U}|\pi)$ for any $g\in G$.
  \item $\overline{D}(G,\mathcal{U}\vee\mathcal{V}|\pi)=\max\{\overline{D}(G,\mathcal{U}|\pi),\overline{D}(G,\mathcal{V}|\pi)\}$.
  \item
  \begin{equation} \label{eq1}
  \begin{aligned}
  &\max\{\underline{D}(G,\mathcal{U}|\pi),\underline{D}(G,\mathcal{V}|\pi)\}\leq \underline{D}(G,\mathcal{U}\vee\mathcal{V}|\pi)\\
  \leq &\max\{\underline{D}(G,\mathcal{U}|\pi),\underline{D}(G,\mathcal{V}|\pi),\min\{\overline{D}(G,\mathcal{U}|\pi),\overline{D}(G,\mathcal{V}|\pi)\}\}.
  \end{aligned}
  \end{equation}
\end{enumerate}
\end{prop}
\begin{proof}
If $\mathcal{U},\mathcal{V}\in \mathcal{C}_X^o$ with $\mathcal{U}\preceq\mathcal{V}$, then $\mathcal{N}(\mathcal{U}|\pi)\leq \mathcal{N}(\mathcal{V}|\pi)$. By the definitions of the relative upper and lower entropy dimension, we get (1).

Since $G$ is a group, we have that $\mathcal{N}(\mathcal{U}|\pi)=\mathcal{N}(g^{-1}\mathcal{U}|\pi)$, for all $\mathcal{U}\in\mathcal{C}_X^o$ and $g\in G$. Then from the definitions of the relative upper and lower entropy dimension, (2) follows.

For (3), we can easily have $\max\{\overline{D}(G,\mathcal{U}|\pi),\overline{D}(G,\mathcal{V}|\pi)\}\leq \overline{D}(G,\mathcal{U}\vee\mathcal{V}|\pi)$. It is sufficient to show $\overline{D}(G,\mathcal{U}\vee\mathcal{V}|\pi)\leq \max\{\overline{D}(G,\mathcal{U}|\pi),\overline{D}(G,\mathcal{V}|\pi)\}$. If $\alpha>\max\{\overline{D}(G,\mathcal{U}|\pi),\overline{D}(G,\mathcal{V}|\pi)\}$, then
\[
\limsup\limits_{n\rightarrow \infty}\frac{\log \mathcal{N}(\bigvee_{g\in F_n}g^{-1}\mathcal{U}|\pi)}{|F_n|^\alpha}=0\mbox{ and }\limsup\limits_{n\rightarrow \infty}\frac{\log \mathcal{N}(\bigvee_{g\in F_n}g^{-1}\mathcal{V}|\pi)}{|F_n|^\alpha}=0.
\]
Since $\mathcal{N}(\bigvee_{g\in F_n}g^{-1}(\mathcal{U}\vee \mathcal{V})|\pi)\leq \mathcal{N}(\bigvee_{g\in F_n}g^{-1}\mathcal{U}|\pi)\cdot \mathcal{N}(\bigvee_{g\in F_n}g^{-1}\mathcal{V}|\pi)$, we have
\[
\begin{aligned}
\overline{h}(G,\mathcal{U}\vee\mathcal{V},\alpha|\pi)&=\limsup\limits_{n\rightarrow \infty}\frac{\log \mathcal{N}(\bigvee_{g\in F_n}g^{-1}(\mathcal{U}\vee \mathcal{V})|\pi)}{|F_n|^\alpha}\\
&\leq \limsup\limits_{n\rightarrow \infty}\frac{\log \mathcal{N}(\bigvee_{g\in F_n}g^{-1}\mathcal{U}|\pi)+\log \mathcal{N}(\bigvee_{g\in F_n}g^{-1}\mathcal{V}|\pi)}{|F_n|^\alpha}\\
&\leq \limsup\limits_{n\rightarrow \infty}\frac{\log \mathcal{N}(\bigvee_{g\in F_n}g^{-1}\mathcal{U}|\pi)}{|F_n|^\alpha}+\limsup\limits_{n\rightarrow \infty}\frac{\log \mathcal{N}(\bigvee_{g\in F_n}g^{-1}\mathcal{V}|\pi)}{|F_n|^\alpha}=0
\end{aligned}
\]
This implies $\overline{D}(G,\mathcal{U}\vee\mathcal{V}|\pi)\leq\alpha$. Since $\alpha$ is arbitrary, we have
\[
\overline{D}(G,\mathcal{U}\vee\mathcal{V}|\pi)\leq \max\{\overline{D}(G,\mathcal{U}|\pi),\overline{D}(G,\mathcal{V}|\pi)\}.
\]
For (4), the first inequality is obvious from (1). We now show \[\underline{D}(G,\mathcal{U}\vee\mathcal{V}|\pi)\leq \max\{\underline{D}(G,\mathcal{U}|\pi),\underline{D}(G,\mathcal{V}|\pi),\min\{\overline{D}(G,\mathcal{U}|\pi),\overline{D}(G,\mathcal{V}|\pi)\}\}.\]
If $\alpha>\max\{\underline{D}(G,\mathcal{U}|\pi),\underline{D}(G,\mathcal{V}|\pi),\min\{\overline{D}(G,\mathcal{U}|\pi),\overline{D}(G,\mathcal{V}|\pi)\}\}$, without loss of generality, we assume $\alpha>\overline{D}(G,\mathcal{V}|\pi)$ and $\alpha>\underline{D}(G,\mathcal{U}|\pi)$. Then
\[
\limsup\limits_{n\rightarrow \infty}\frac{\log \mathcal{N}(\bigvee_{g\in F_n}g^{-1}\mathcal{V}|\pi)}{|F_n|^\alpha}=0\mbox{ and }\liminf\limits_{n\rightarrow \infty}\frac{\log \mathcal{N}(\bigvee_{g\in F_n}g^{-1}\mathcal{U}|\pi)}{|F_n|^\alpha}=0.
\]
Hence
\[
\begin{aligned}
\underline{h}(G,\mathcal{U}\vee\mathcal{V},\alpha|\pi)&=\liminf\limits_{n\rightarrow \infty}\frac{\log \mathcal{N}(\bigvee_{g\in F_n}g^{-1}(\mathcal{U}\vee \mathcal{V})|\pi)}{|F_n|^\alpha}\\
&\leq \liminf\limits_{n\rightarrow \infty}\frac{\log \mathcal{N}(\bigvee_{g\in F_n}g^{-1}\mathcal{U}|\pi)+\log \mathcal{N}(\bigvee_{g\in F_n}g^{-1}\mathcal{V}|\pi)}{|F_n|^\alpha}\\
&\leq \liminf\limits_{n\rightarrow \infty}\frac{\log \mathcal{N}(\bigvee_{g\in F_n}g^{-1}\mathcal{U}|\pi)}{|F_n|^\alpha}+\limsup\limits_{n\rightarrow \infty}\frac{\log \mathcal{N}(\bigvee_{g\in F_n}g^{-1}\mathcal{V}|\pi)}{|F_n|^\alpha}=0
\end{aligned}
\]
This implies $\underline{D}(G,\mathcal{U}\vee\mathcal{V}|\pi)<\alpha$. Since $\alpha$ is arbitrary, we have
\[\underline{D}(G,\mathcal{U}\vee\mathcal{V}|\pi)\leq \max\{\underline{D}(G,\mathcal{U}|\pi),\underline{D}(G,\mathcal{V}|\pi),\min\{\overline{D}(G,\mathcal{U}|\pi),\overline{D}(G,\mathcal{V}|\pi)\}\}.\]
\end{proof}
From Proposition \ref{p1}, we can get the following result.
\begin{prop}\label{p2}
Let $\pi:(X,G)\rightarrow (Y,G)$ be a factor map between $G$-systems. If $\{\mathcal{U}_n\}$ is sequence chosen from $\mathcal{C}_X^o$ with $\lim\limits_{n\rightarrow \infty}diam(\mathcal{U}_n)=0$, then
\[
\lim\limits_{n\rightarrow \infty}\overline{D}(G,\mathcal{U}_n|\pi)=\overline{D}(X,G|\pi)\mbox{ and } \lim\limits_{n\rightarrow \infty}\underline{D}(G,\mathcal{U}_n|\pi)=\underline{D}(X,G|\pi).
\]
In particular, if $G$ is an abelian group and $\mathcal{U}$ is a generating open cover of $X$ (i.e. $\lim\limits_{n\rightarrow \infty}diam(\bigvee\limits_{g\in F_n}g^{-1}\mathcal{U})=0$), then
\[
\overline{D}(G,\mathcal{U}|\pi)=\overline{D}(X,G|\pi).
\]
\end{prop}
\begin{proof}
We only need to consider the upper case. Since $\lim\limits_{n\rightarrow \infty}diam(\mathcal{U}_n)=0$, by (1) of Proposition \ref{p1}, we have \[\lim\limits_{n\rightarrow \infty}\overline{D}(G,\mathcal{U}_n|\pi)=\sup\limits_n\overline{D}(G,\mathcal{U}_n|\pi).\] Therefore, \[\overline{D}(X,G|\pi)=\sup\limits_n\overline{D}(G,\mathcal{U}_n|\pi)=\lim\limits_{n\rightarrow \infty}\overline{D}(G,\mathcal{U}_n|\pi).\]
For the particular case, $\lim\limits_{n\rightarrow \infty}diam(\bigvee\limits_{g\in F_n}g^{-1}\mathcal{U})=0$. By (2) and (3) of Proposition \ref{p1}, we have
\[\overline{D}(G,\bigvee\limits_{g\in F_n}g^{-1}\mathcal{U}|\pi)=\overline{D}(G,\mathcal{U}|\pi).\] Therefore, \[\overline{D}(X,G|\pi)=\sup\limits_n\overline{D}(G,\bigvee\limits_{g\in F_n}g^{-1}\mathcal{U}|\pi)=\sup\limits_n\overline{D}(G,\mathcal{U}|\pi)=\overline{D}(G,\mathcal{U}|\pi).\]
\end{proof}

Let $(X,G)$ be a $G$-system. A cover $\{U,V\}$ of $X$ which consists of two non-dense open sets of $X$ is called \textit{a standard cover} of $X$. Denote by $\mathcal{C}_X^s$ the set of all standard covers of $X$. The following proposition shows the relative upper entropy dimension with respect to standard covers determine the relative upper entropy dimension of the extension.
\begin{prop}\label{p3}
Let $\pi:(X,G)\rightarrow (Y,G)$ be a factor map between $G$-systems. Then
\[
\overline{D}(X,G|\pi)=\sup\{\overline{D}(G,\mathcal{W}|\pi):\mathcal{W}\in\mathcal{C}_X^s\}.
\]
\end{prop}
\begin{proof}
We follow the argument in the proof of \cite[Proposition 1]{Bl}. Since $\mathcal{C}_X^o\supseteq \mathcal{C}_X^s$, we have
\[
\overline{D}(X,G|\pi)\geq\sup\{\overline{D}(G,\mathcal{W}|\pi):\mathcal{W}\in\mathcal{C}_X^s\}.
\]
We now show the opposite inequality. Let $\mathcal{U}=\{U_1,U_2,\cdots,U_k\}\in \mathcal{C}_X^o$. By the compactness of $X$, there exists a finite open subcover $\{W_1,W_2,\cdots,W_p\}$ such that for any $1\leq j\leq p$, there exists $1\leq i\leq k$ such that $\overline{W_j}\subseteq U_i$. Let $F_j=\overline{W_j}$, $j=1,2,\cdots,p$ and $V_i=\bigcap\limits_{j\in\{j:F_j\subseteq U_i\}}F_j^c$. Since $U_i\supseteq V_i^c$, $\mathcal{U}_i=\{U_i,V_i\}$ is an open cover of $X$ and $\bigcap\limits_{i=1}^kV_i=\emptyset$. Since every element of $\bigvee_{i=1}^k\mathcal{U}_i$ is either in some $U_i$ or equal to $\bigcap\limits_{i=1}^kV_i=\emptyset$, we have $\bigvee_{i=1}^k\mathcal{U}_i\succeq\mathcal{U}$. Note that
\[
\overline{D}(G,\mathcal{U}|\pi)\leq \overline{D}(G,\bigvee_{i=1}^k\mathcal{U}_i|\pi)=\max\{\overline{D}(G,\mathcal{U}_i|\pi):i=1,2,\cdots,k\}.
\]
hence there exists $r\in\{1,2,\cdots,k\}$ such that $\overline{D}(G,\mathcal{U}_r|\pi)\geq \overline{D}(G,\mathcal{U}|\pi)$.

Now it is easy to shrink $\mathcal{U}_k=\{U_r,V_r\}$ into a standard cover $\mathcal{W}$ with $\overline{D}(G,\mathcal{W}|\pi)\geq \overline{D}(G,\mathcal{U}|\pi)$. First if $U_r$ is dense, then for $\varepsilon>0$, there exists $x\in U_r$ such that \[F=\{y:d(x,y)\leq\varepsilon\}\subseteq V_r.\]
Let $W_1=U_r\setminus F$. Since $F$ has nonempty interior, the set $W_1$ is not dense and $\{W_1,V_r\}$ is an open cover of $X$ as $F\subseteq V_r$. Since $\{W_1,V_r\}\succeq\mathcal{U}_r$, we have
\[
\overline{D}(G,\{W_1,V_r\}|\pi)\geq\overline{D}(G,\mathcal{U}_r|\pi).
\]
By doing the same operation with $V_r$ if $V_r$ is dense, one obtains a standard cover $\mathcal{W}=\{W_1,W_2\}$ such that
\[
\overline{D}(G,\mathcal{W}|\pi)\geq\overline{D}(G,\mathcal{U}|\pi).
\]
Since $\mathcal{U}$ is arbitrary, $\sup\{\overline{D}(G,\mathcal{W}|\pi):\mathcal{W}\in\mathcal{C}_X^s\}\geq \overline{D}(X,G|\pi)$.
\end{proof}
\section{Relative entropy dimensions via entropy generating sets}\label{sec4}
Given a strict increasing F${\o}$lner sequence $\{F_n\}$ of $G$, we suppose that a set $S\subseteq G$ satisfies that
\begin{equation}\label{4.1}
\lim\limits_{n\rightarrow\infty}|S\cap F_n|=+\infty.
\end{equation}
We denote by $\mathcal{I}(G)$ the set of the subset $S$ of $G$ satisfying with (\ref{4.1}). For $\alpha\geq0$ and $S\in \mathcal{I}(G)$, we define
\[
\overline{D}(S,\alpha)=\limsup\limits_{n\rightarrow \infty} \frac{|S\cap F_n|}{|F_n|^\alpha}\mbox{ and } \underline{D}(S,\alpha)=\liminf\limits_{n\rightarrow \infty} \frac{|S\cap F_n|}{|F_n|^\alpha}.
\]
It is clear that $\overline{D}(S,\alpha)$ does not decrease as $\alpha$ decreases and $\overline{D}(S,\alpha)\notin\{0,+\infty\}$ for at most one $\alpha\geq 0$. We define the upper dimension of $S$ by
\[
\overline{D}(S)=\inf\{\alpha\geq0:\overline{D}(S,\alpha)=0\}=\sup\{\alpha\geq0:\overline{D}(S,\alpha)=\infty\}.
\]
Similarly, we define the lower dimension of $S$ by
\[
\underline{D}(S)=\inf\{\alpha\geq0:\underline{D}(S,\alpha)=0\}=\sup\{\alpha\geq0:\underline{D}(S,\alpha)=\infty\}.
\]
Clearly $0\leq\underline{D}(S)\leq\overline{D}(S)\leq1$. When $\overline{D}(S)=\underline{D}(S)=\alpha$, we say the set $S$ has dimension $\alpha$, denoted by $D(S)=\alpha$. If $S$ has a positive density with respect to $\{F_n\}$ i.e. $\lim\limits_{n\rightarrow \infty}\frac{|F_n\cap S|}{|F_n|}>0$, then $D(S)=1$.

In the following, we will investigate the dimension of a special kind of subsets of $G$, which we call the relative entropy generating sets.

Let $\pi:(X,G)\rightarrow (Y,G)$ be a factor map between $G$-systems and $\mathcal{U}\in\mathcal{C}_X^o$. We say $S\in \mathcal{I}(G)$ is a relative entropy generating set of $\mathcal{U}$ relevant to $\pi$ if
\[
\liminf\limits_{n\rightarrow \infty}\frac{1}{|F_n\cap S|}\log\mathcal{N}\left(\bigvee\limits_{g\in F_n\cap S}g^{-1}\mathcal{U}|\pi\right)>0.
\]
Denote by $\mathcal{E}(G,\mathcal{U}|\pi)$ the set of all relative entropy generating sets of $\mathcal{U}$ relevant to $\pi$ and by $\mathcal{P}(G,\mathcal{U}|\pi)$ the set of $S\in \mathcal{I}(G)$ with the property that
\[
\limsup\limits_{n\rightarrow \infty}\frac{1}{|F_n\cap S|}\log\mathcal{N}\left(\bigvee\limits_{g\in F_n\cap S}g^{-1}\mathcal{U}|\pi\right)>0.
\]
In other words, $\mathcal{P}(G,\mathcal{U}|\pi)$ is the set of subsets of $G$ along which $\mathcal{U}$ has positive relative upper entropy.
\begin{defn}
Let $\pi:(X,G)\rightarrow (Y,G)$ be a factor map between $G$-systems and $\mathcal{U}\in \mathcal{C}_X^o$. We define
\begin{equation*}
\overline{D}_e(G,\mathcal{U}|\pi)=
\begin{cases}
\sup\limits_{S\in \mathcal{E}(G,\mathcal{U}|\pi)}\overline{D}(S), &\mbox{ if }\mathcal{E}(G,\mathcal{U}|\pi)\neq \emptyset,\\
0, &\mbox{ if }\mathcal{E}(G,\mathcal{U}|\pi)= \emptyset,
\end{cases}
\end{equation*}
\begin{equation*}
\overline{D}_p(G,\mathcal{U}|\pi)=
\begin{cases}
\sup\limits_{S\in \mathcal{P}(G,\mathcal{U}|\pi)}\overline{D}(S), & \mbox{ if }\mathcal{P}(G,\mathcal{U}|\pi)\neq \emptyset,\\
0, & \mbox{ if }\mathcal{P}(G,\mathcal{U}|\pi)= \emptyset.
\end{cases}
\end{equation*}
\end{defn}
Similarly, we can define $\underline{D}_e(G,\mathcal{U}|\pi)$ and $\underline{D}_p(G,\mathcal{U}|\pi)$ by changing the upper dimension into the lower dimension.
\begin{defn}
Let $\pi:(X,G)\rightarrow (Y,G)$ be a factor map between $G$-systems. We define
\[
\overline{D}_e(X,G|\pi)=\sup\limits_{\mathcal{U}\in\mathcal{C}_X^o}\overline{D}_e(G,\mathcal{U}|\pi), \overline{D}_p(X,G|\pi)=\sup\limits_{\mathcal{U}\in\mathcal{C}_X^o}\overline{D}_p(G,\mathcal{U}|\pi).
\]
\end{defn}
Similarly, we can define $\underline{D}_e(X,G|\pi)$ and $\underline{D}_p(X,G|\pi)$.

The following proposition explains why we define the relative entropy generating set as $\liminf$ instead of $\limsup$.
\begin{prop}
Let $\pi:(X,G)\rightarrow (Y,G)$ be a factor map between $G$-systems and $\mathcal{U}\in \mathcal{C}_X^o$. Then
\begin{equation*}
\overline{D}_p(G,\mathcal{U}|\pi)=
\begin{cases}
1, & \mbox{ if }\mathcal{P}(G,\mathcal{U}|\pi)\neq \emptyset,\\
0, & \mbox{ if }\mathcal{P}(G,\mathcal{U}|\pi)= \emptyset.
\end{cases}
\end{equation*}
\end{prop}
\begin{proof}
We assume $\mathcal{P}(G,\mathcal{U}|\pi)\neq\emptyset$. Then there are $a>0$ and $S\in \mathcal{I}(G)$ such that
\[
\limsup\limits_{n\rightarrow\infty}\frac{1}{|F_n\cap S|}\log\mathcal{N}\left(\bigvee\limits_{g\in F_n\cap S}g^{-1}\mathcal{U}|\pi\right)=a.
\]
Next we take $1\leq n_1<n_2<\cdots$ such that \[2|F_{n_j}\cap S|\leq|F_{n_{j+1}}|\leq|F_{n_{j+2}}\cap S| \mbox{ for all } j\in\mathbb{N}
\]
and
\[
\limsup\limits_{j\rightarrow\infty}\frac{1}{|F_{n_j}\cap S|}\log\mathcal{N}\left(\bigvee\limits_{g\in F_{n_j}\cap S}g^{-1}\mathcal{U}|\pi\right)=a.
\]
Put
\[
F=F_{n_1}\cup S\cup\bigcup\limits_{i=1}^\infty F_{n_{i+1}}\setminus(F_{n_i}\cap S).
\]
Then $F\in \mathcal{I}(G)$ and
\[
\begin{aligned}
\limsup\limits_{n\rightarrow \infty}\frac{1}{|F\cap F_n|}\log\mathcal{N}\left(\bigvee_{g\in F\cap F_n}g^{-1}\mathcal{U}|\pi\right)&\geq\limsup\limits_{j\rightarrow \infty}\frac{1}{|F\cap((S\cap F_{n_j})\cup F_{n_{j-1}})|}\log\mathcal{N}\left(\bigvee\limits_{g\in F_{n_j}\cap S}g^{-1}\mathcal{U}|\pi\right)\\
&\geq \limsup\limits_{j\rightarrow \infty}\frac{1}{|(S\cap F_{n_j})\cup F_{n_{j-1}}|}\log\mathcal{N}\left(\bigvee\limits_{g\in F_{n_j}\cap S}g^{-1}\mathcal{U}|\pi\right)\\
&\geq \limsup\limits_{j\rightarrow \infty}\frac{1}{2|S\cap F_{n_j}|}\log\mathcal{N}\left(\bigvee\limits_{g\in F_{n_j}\cap S}g^{-1}\mathcal{U}|\pi\right)=\frac{a}{2}>0,
\end{aligned}
\]
therefore, $F\in \mathcal{P}(G,\mathcal{U}|\pi)$. Since $|F_{n_{j+1}}|\geq2|F_{n_j}\cap S|$ for each $j\in\mathbb{N}$, it is easy to see that the upper density of $F$ is $\frac{1}{2}$, hence $\overline{D}(F)=1$. This implies $\overline{D}_p(G,\mathcal{U}|\pi)=1$.
\end{proof}
\section{Some relationships among relative entropy dimensions}\label{sec5}
In the following part, we investigate the interrelations among the relative entropy dimensions appeared in Section \ref{sec3} and \ref{sec4}.
\begin{prop}\label{p4}
Let $\pi:(X,G)\rightarrow (Y,G)$ be a factor map between $G$-systems and $\mathcal{U}\in \mathcal{C}_X^o$. Then
\[
\underline{D}_e(G,\mathcal{U}|\pi)\leq\overline{D}_e(G,\mathcal{U}|\pi)\leq\underline{D}_p(G,\mathcal{U}|\pi)\leq\overline{D}(G,\mathcal{U}|\pi).
\]
\end{prop}
\begin{proof}
\begin{enumerate}[1)]
  \item $\underline{D}_e(G,\mathcal{U}|\pi)\leq\overline{D}_e(G,\mathcal{U}|\pi)$ is obvious by the definition.
  \item To show $\overline{D}_e(G,\mathcal{U}|\pi)\leq\underline{D}_p(G,\mathcal{U}|\pi)$, it is sufficient to assume that $\overline{D}_e(G,\mathcal{U}|\pi)>0$. Let $\alpha\in (0,\overline{D}_e(G,\mathcal{U}|\pi))$ be given. Then there exists $S\in\mathcal{E}(G,\mathcal{U}|\pi)$ such that $\overline{D}(S)>\alpha$, that is,
      \[
      \limsup\limits_{n\rightarrow+\infty}\frac{|S\cap F_n|}{|F_n|^\alpha}=+\infty.
      \]
  Therefore, we have that
  \begin{equation}\label{eq2}
   \limsup\limits_{n\rightarrow+\infty}\frac{|S\cap F_n|}{|S\cap F_n|+|F_n|^\alpha}=1.
  \end{equation}
  For any $i\in\mathbb{N}$, choose $M_i\subseteq F_i\setminus F_{i-1}$ such that $|M_i|=\lfloor|F_i|^\alpha\rfloor-\lfloor|F_{i-1}|^\alpha\rfloor$, where $\lfloor x\rfloor$ denotes the maximal integer not greater than $x$. Let $F=S\cup \bigcup\limits_{i=1}^\infty M_i$. Then $F\in\mathcal{I}(G)$ and it is easy to see that $\underline{D}(F)\geq \alpha$. Write $F=\{t_1,t_2,\cdots\}$ and order it as $M_1\cup(F_1\cap S),M_2\cup((F_2\setminus F_1)\cap S),\cdots,M_n\cup((F_n\setminus F_{n-1})\cap S), \cdots, F\setminus(\bigcup\limits_{i=1}^\infty(M_n\cup(F_n\cap S)))$. Then for any large enough $n\in\mathbb{N}$, there exists an unique $m(n)\in\mathbb{N}$ such that $m(n)=\max\{i:t_i\in F\cap F_n\}$. Since
  \[
  F_n\cap S\subseteq \{t_1,t_2,\cdots,t_{m(n)}\}\subseteq (F_n\cap S)\overset{n}{\underset{i=1}{\cup}}M_i,
  \]
  we have that $|S\cap F_n|\leq m(n)\leq |S\cap F_n|+|F_n|^\alpha$, hence by (\ref{eq2}), we have that
  \begin{equation}\label{eq}
  \limsup\limits_{n\rightarrow\infty}\frac{|S\cap F_n|}{m(n)}=1.
  \end{equation}
  Since $S\in \mathcal{E}(G,\mathcal{U}|\pi)$, we have that
  \[
  \begin{aligned}
  \limsup\limits_{n\rightarrow +\infty}\frac{1}{|F\cap F_n|}\log\mathcal{N}\left(\bigvee\limits_{g\in F\cap F_n}g^{-1}\mathcal{U}|\pi\right)\geq&\limsup\limits_{n\rightarrow +\infty}\frac{1}{m(n)}\log\mathcal{N}\left(\bigvee\limits_{i=1}^{m(n)}t_i^{-1}\mathcal{U}|\pi\right)\\
  \geq&\limsup\limits_{n\rightarrow +\infty}\frac{\log\mathcal{N}\left(\bigvee\limits_{g\in F_n\cap S}g^{-1}\mathcal{U}|\pi\right)}{|S\cap F_n|}\cdot\frac{|S\cap F_n|}{m(n)}\\
  \geq &\liminf\limits_{n\rightarrow +\infty}\frac{\log\mathcal{N}\left(\bigvee\limits_{g\in F_n\cap S}g^{-1}\mathcal{U}|\pi\right)}{|S\cap F_n|}\cdot\limsup\limits_{n\rightarrow +\infty}\frac{|S\cap F_n|}{m(n)}\\
  \xlongequal{(\ref{eq})}&\liminf\limits_{n\rightarrow +\infty}\frac{\log\mathcal{N}\left(\bigvee\limits_{g\in F_n\cap S}g^{-1}\mathcal{U}|\pi\right)}{|S\cap F_n|}>0.
  \end{aligned}
  \]
  This implies $F\in \mathcal{P}(G,\mathcal{U}|\pi)$. Hence $\underline{D}_p(G,\mathcal{U}|\pi)\geq \underline{D}(F)\geq\alpha$. Since $\alpha$ is arbitrary chosen from $\alpha\in (0,\overline{D}_e(G,\mathcal{U}|\pi))$, we have $\overline{D}_e(G,\mathcal{U}|\pi)\leq\underline{D}_p(G,\mathcal{U}|\pi)$.
  \item If $\underline{D}_p(G,\mathcal{U}|\pi)\leq \overline{D}(G,\mathcal{U}|\pi)$ were not true, then there exists $\alpha\in(0,1)$ such that
  \[
  \underline{D}_p(G,\mathcal{U}|\pi)>\alpha>\overline{D}(G,\mathcal{U}|\pi).
  \]
  On the one hand, since $\alpha>\overline{D}(G,\mathcal{U}|\pi)$, we have that
  \begin{equation}\label{eq3}
  \limsup\limits_{n\rightarrow +\infty}\frac{1}{|F_n|^\alpha}\log\mathcal{N}\left(\bigvee\limits_{g\in F_n}g^{-1}\mathcal{U}|\pi\right)=0.
  \end{equation}
  On the other hand, since $\underline{D}_p(G,\mathcal{U}|\pi)>\alpha$, there exists $S\in \mathcal{P}(G,\mathcal{U}|\pi)$ such that $\underline{D}(S)>\alpha$, that is,
  \[
  \liminf\limits_{n\rightarrow \infty}\frac{|F_n\cap S|}{|F_n|^\alpha}=+\infty.
  \]
  Hence there is $c>0$ such that $\frac{|F_n\cap S|}{|F_n|^\alpha}\geq c$ for all sufficiently large $n$. Now
  \[
  \begin{split}
  \limsup\limits_{n\rightarrow +\infty}\frac{1}{|F_n|^\alpha}\log\mathcal{N}\left(\bigvee\limits_{g\in F_n}g^{-1}\mathcal{U}|\pi\right)
  \geq&\limsup\limits_{n\rightarrow +\infty}\frac{1}{|F_n|^\alpha}\log\mathcal{N}\left(\bigvee\limits_{g\in F_n\cap S}g^{-1}\mathcal{U}|\pi\right)\\
  \geq&\limsup\limits_{n\rightarrow +\infty}\frac{\log\mathcal{N}\left(\bigvee\limits_{g\in F_n\cap S}g^{-1}\mathcal{U}|\pi\right)}{|F_n\cap S|}\cdot \frac{|F_n\cap S|}{|F_n|^\alpha}\\
  \geq&\limsup\limits_{n\rightarrow +\infty}\frac{\log\mathcal{N}\left(\bigvee\limits_{g\in F_n\cap S}g^{-1}\mathcal{U}|\pi\right)}{|F_n\cap S|}\cdot c>0,\\
  &(\mbox{ since }S\in\mathcal{P}(G,\mathcal{U}|\pi))
  \end{split}
  \]
  which contradicts (\ref{eq3}).
\end{enumerate}
\end{proof}
\begin{prop}\label{p5}
Let $\pi:(X,G)\rightarrow (Y,G)$ be a factor map between $G$-systems and $\mathcal{U}\in \mathcal{C}_X^o$. Then
\[
\underline{D}_e(G,\mathcal{U}|\pi)\leq\underline{D}(G,\mathcal{U}|\pi)\leq\overline{D}(G,\mathcal{U}|\pi).
\]
\end{prop}
\begin{proof}
By Definition \ref{d1}, it is clear that $\underline{D}(G,\mathcal{U}|\pi)\leq\overline{D}(G,\mathcal{U}|\pi).$ Now we will show that $\underline{D}_e(G,\mathcal{U}|\pi)\leq\underline{D}(G,\mathcal{U}|\pi)$. If $\underline{D}_e(G,\mathcal{U}|\pi)=0$, then the inequality holds obviously, otherwise, let $\alpha\in(0,\underline{D}_e(G,\mathcal{U}|\pi))$. Then there exists $S\in\mathcal{E}(G,\mathcal{U}|\pi)$ such that $\underline{D}(S)>\alpha$, that is,
\[
\liminf\limits_{n\rightarrow\infty}\frac{|F_n\cap S|}{|F_n|^\alpha}=+\infty.
\]
Hence there is $c>0$ such that $|S\cap F_n|\geq c|F_n|^\alpha$ for sufficiently large $n$. Since $S\in\mathcal{E}(G,\mathcal{U}|\pi)$, we have
\[
\liminf\limits_{n\rightarrow +\infty}\frac{1}{|F_n|^\alpha}\log\mathcal{N}\left(\bigvee\limits_{g\in F_n}g^{-1}\mathcal{U}|\pi\right)\geq c\liminf\limits_{n\rightarrow +\infty}\frac{1}{|F_n\cap S|}\log\mathcal{N}\left(\bigvee\limits_{g\in F_n\cap S}g^{-1}\mathcal{U}|\pi\right)>0.
\]
This implies $\underline{D}(G,\mathcal{U}|\pi)\geq\alpha$. Finally, since $\alpha\in(0,\underline{D}_e(G,\mathcal{U}|\pi))$ is arbitrary, we get $\underline{D}_e(G,\mathcal{U}|\pi)\leq\underline{D}(G,\mathcal{U}|\pi)$.
\end{proof}
Let $k\geq 2$ and $\pi:(X,G)\rightarrow (Y,G)$ be a factor map between $G$-systems. Let $A_1,A_2,\cdots,A_k$ be $k$ subsets of $X$ and $W\subseteq G$. We say $\{A_1,A_2,\cdots,A_k\}$ is \textit{independent} along $W$ relevant to $\pi$, if there is $y\in Y$ such that for any $s\in\{1,2,\cdots,k\}^W$, we have
$\bigcap\limits_{w\in W} w^{-1}A_{s(w)}\cap \pi^{-1}(y)\neq\emptyset$ (see \cite{KL1}).
\begin{lem}\label{l1}(Cf. \cite{Zhou} for $G=\mathbb{Z}$)
Let $\pi:(X,G)\rightarrow (Y,G)$ be a factor map between $G$-systems, and let $A_1,A_2$ be two disjoint non-empty closed subsets of $X$, $\mathcal{U}=\{A_1^c,A_2^c\}$. For any $\alpha\in (0,1]$, $0<\eta<\alpha$ and $c>0$, there exists $N\in\mathbb{N}$ (depending on $\alpha,\eta,c$) such that if a finite subset $B$ of $G$ with $|B|\geq N$ satisfies \[\mathcal{N}\left(\bigvee\limits_{g\in B}g^{-1}\mathcal{U}|\pi\right)>e^{c|B|^\alpha},\] then there exists $W\subseteq B$ with $|W|\geq|B|^\eta$ such that $\{A_1,A_2\}$ is independent along $W$ relevant to $\pi$.
\end{lem}
\begin{proof}
The proof of the lemma is similar to the proof of \cite[Lemma 3.7]{DHP}.
\end{proof}
With the help of the above lemma, we have the following theorem.
\begin{thm}\label{t1}
Let $\pi:(X,G)\rightarrow (Y,G)$ be a factor map between $G$-systems, and let $A_1,A_2$ be two disjoint non-empty closed subsets of $X$, $\mathcal{U}=\{A_1^c,A_2^c\}$. Then there exists a set $F\in \mathcal{E}(G,\mathcal{U}|\pi)$ such that $\overline{D}(F)=\overline{D}(G,\mathcal{U}|\pi)$, when $\mathcal{E}(G,\mathcal{U}|\pi)$ is non-empty. Moreover, $\overline{D}(G,\mathcal{U}|\pi)=\overline{D}_e(G,\mathcal{U}|\pi)$.
\end{thm}
\begin{proof}
If $\overline{D}(G,\mathcal{U}|\pi)=0$, then by Proposition \ref{p4}, any set in $\mathcal{E}(G,\mathcal{U}|\pi)$ will have upper dimension zero.

Now we assume that $\overline{D}(G,\mathcal{U}|\pi)>0$ and let $\{\alpha_j\}\subseteq (0,\overline{D}(G,\mathcal{U}|\pi))$ be a sequence of strictly increasing real number such that $\lim\limits_{j\rightarrow \infty}\alpha_j=\overline{D}(G,\mathcal{U}|\pi)$. Then we can choose $a>0$ so that
\[
\limsup\limits_{n\rightarrow \infty}\frac{1}{|F_n|^{\alpha_j}}\log\mathcal{N}(\bigvee\limits_{g\in F_n}g^{-1}\mathcal{U}|\pi)>a, \mbox{ for all }j\in\mathbb{N}.
\]
Let $\alpha_{j-1}<\eta_j<\alpha_j$ for each $j\in\mathbb{N}$. Then by Lemma \ref{l1}, if there exists $N_j\in \mathbb{N}$ such that for every finite set $C\subseteq G$ with $|C|\geq N_j$ and $\mathcal{N}\left(\bigvee\limits_{g\in C}g^{-1}\mathcal{U}|\pi\right)\geq e^{\frac{a}{2}|C|^{\alpha_j}}$, then we can find $W\subseteq C$ with $|W|\geq|C|^{\eta_j}$ such that $\{A_1,A_2\}$ is independent along $W$ relevant to $\pi$.

Take $1=n_1<n_2<\cdots$ such that
\[
|F_{n_{j+1}}\setminus F_{n_j}|^{\eta_j}\geq j|F_{n_j}|+N_j\mbox{ and }\mathcal{N}\left(\bigvee\limits_{g\in F_{n_{j+1}}\setminus F_{n_j}}g^{-1}\mathcal{U}|\pi\right)\geq e^{\frac{a}{2}|F_{n_{j+1}}\setminus F_{n_j}|^{\alpha_j}}.
\]
Now for each $j\in\mathbb{N}$, there exists $W_j\subseteq F_{n_{j+1}}\setminus F_{n_j}$ such that $|W_j|\geq|F_{n_{j+1}}\setminus F_{n_j}|^{\eta_j}$ and $\{A_1,A_2\}$ is independent along $W_j$ relevant to $\pi$, that is, there exists $y_j\in Y$ such that $\forall s\in \{1,2\}^{W_j}$, we have
\[
\bigcap\limits_{w\in W_j} w^{-1}A_{s(w)}\cap \pi^{-1}(y_j)\neq\emptyset.
\]

For any non-empty set $B\subseteq W_j$ and $s=(s(z))_{z\in B}\in\{1,2\}^B$, we can find $x_s\in \bigcap\limits_{z\in B}z^{-1}A_{s(z)}\cap\pi^{-1}(y_j)$. Let $X_B=\{x_s:s\in\{1,2\}^B\}$. It is clear that for any $l\in\{1,2\}^B$, we have
\[
\left|\bigcap\limits_{z\in B}z^{-1}A^c_{s(z)}\cap\pi^{-1}(y_j)\cap X_B\right|\leq (2-1)^{|B|}=1.
\]
Combining this fact with $|X_B|=2^{|B|}$, we get
\begin{equation}\label{eq4}
\mathcal{N}\left(\bigvee\limits_{z\in B}z^{-1}\mathcal{U},\pi^{-1}(y_j)\right)\geq 2^{|B|},\mbox{ for all } B\subseteq W_j.
\end{equation}

Put $F=\bigcup\limits_{i=1}^\infty W_j$ and order it as $W_1W_2\cdots$ (for example, the $|W_1|$ elements of $F$ is the whole $W_1$). Then $F\in\mathcal{I}(G)$ and for $n\in\mathbb{N}$ with $|F_n\cap F|\geq|W_1|$, there exists an unique $k(n)\in\mathbb{N}$ such that
\[
\sum\limits_{i=1}^{k(n)} |W_i|\leq |F_n\cap F|<\sum\limits_{i=1}^{k(n)+1} |W_i|.
\]
Now
\[
\begin{split}
\mathcal{N}\left(\bigvee\limits_{g\in F\cap F_n}g^{-1}\mathcal{U}|\pi\right)&\geq\max\left\{\mathcal{N}\left(\bigvee\limits_{w\in W_{k(n)}}w^{-1}\mathcal{U}|\pi\right),\mathcal{N}\left(\bigvee\limits_{w\in W_{k(n)+1}\cap F\cap F_n}w^{-1}\mathcal{U}|\pi\right)\right\}\\
&\overset{(\ref{eq4})}{\geq}\max\{2^{|W_{k(n)}|},2^{|F_n\cap F|-\sum\limits_{i=1}^{k(n)}|W_i|}\}\geq2^{\frac{|F_n\cap F|-\sum\limits_{i=1}^{k(n)-1}|W_i|}{2}}.
\end{split}
\]
Hence
\[
\begin{split}
&\liminf\limits_{n\rightarrow\infty}\frac{1}{|F_n\cap F|}\log\mathcal{N}\left(\bigvee\limits_{g\in F\cap F_n}g^{-1}\mathcal{U}|\pi\right)\geq\liminf\limits_{n\rightarrow\infty}\frac{|F_n\cap F|-\sum\limits_{j=1}^{k(n)-1}|W_j|}{2|F_n\cap F|}\log2\\
\geq &\liminf\limits_{n\rightarrow \infty}\left(\frac{1}{2}-\frac{\sum\limits_{j=1}^{k(n)-1}|F_{n_{j-1}}\setminus F_{n_j}|}{2|W_{k(n)}|}\right)\cdot\log2
\geq \liminf\limits_{n\rightarrow \infty}\left(\frac{1}{2}-\frac{|F_{n_{k(n))}}|}{2|F_{n_{k(n)+1}}\setminus F_{n_{k(n)}}|^{\eta_{k(n)}}}\right)\cdot\log2\\
\geq &\liminf\limits_{n\rightarrow \infty}\left(\frac{1}{2}-\frac{1}{2k(n)}\right)\cdot\log2=\frac{1}{2}\log2>0.
\end{split}
\]
This shows $F\in\mathcal{E}(G,\mathcal{U}|\pi)$.

Note that
\[
\limsup\limits_{n\rightarrow \infty}\frac{|F_n\cap F|}{|F_n|^{\eta_j}}\geq\limsup\limits_{n\rightarrow \infty}\frac{|W_{k(n)}|}{|F_{n_{k(n)}}|^{\eta_j}}\geq\limsup\limits_{n\rightarrow \infty}\frac{|F_{n_{k(n)+1}}\setminus F_{n_{k(n)}}|^{\eta_{k(n)}}}{|F_{n_{k(n)}}|^{\eta_j}}\geq1,
\]
we have $\overline{D}(F)\geq \eta_j$. Hence $\overline{D}(F)=\overline{D}(G,\mathcal{U}|\pi)$. By Proposition \ref{p4}, $\overline{D}(G,\mathcal{U}|\pi)=\overline{D}_e(G,\mathcal{U}|\pi)$. This finishes the proof.
\end{proof}
\begin{thm}
Let $\pi:(X,G)\rightarrow (Y,G)$ be a factor map between $G$-systems. Then
\begin{enumerate}[(1)]
  \item $\overline{D}_e(X,G|\pi)=\underline{D}_p(X,G|\pi)=\overline{D}(X,G|\pi)$.
  \item $\underline{D}_e(X,G|\pi)\leq\underline{D}(X,G|\pi)\leq\overline{D}(X,G|\pi)\leq\overline{D}_p(X,G|\pi)$.
\end{enumerate}
\end{thm}
\begin{proof}
\begin{enumerate}[(1)]
  \item By Proposition \ref{p4}, we have $\overline{D}_e(X,G|\pi)\leq\underline{D}_p(X,G|\pi)\leq\overline{D}(X,G|\pi)$. Now it is sufficient to show that $\overline{D}(X,G|\pi)\leq\overline{D}_e(X,G|\pi)$. By Theorem \ref{t1}, we get that $\overline{D}_e(G,\mathcal{W}|\pi)=\overline{D}(G,\mathcal{W}|\pi)$, for any $\mathcal{W}\in\mathcal{C}_X^s$. Then by Proposition \ref{p3} , we have
\[
\overline{D}(X,G|\pi)=\sup\{\overline{D}(G,\mathcal{W}|\pi):\mathcal{W}\in\mathcal{C}_X^s\}=\sup\{\overline{D}_e(G,\mathcal{W}|\pi):\mathcal{W}\in\mathcal{C}_X^s\}\leq\overline{D}_e(X,G|\pi).
\]
  \item Let $\mathcal{U}\in \mathcal{C}_X^o$. By Definition \ref{d1}, it is clear that $\underline{D}(G,\mathcal{U}|\pi)\leq \overline{D}(G,\mathcal{U}|\pi)$. By Proposition \ref{p5}, we have $\overline{D}_e(G,\mathcal{U}|\pi)\leq\underline{D}(G,\mathcal{U}|\pi)$. Then by (1), $\overline{D}(X,G|\pi)=\underline{D}_p(X,G|\pi)\leq\overline{D}_p(X,G|\pi)$. This finishes the proof of (2).
\end{enumerate}
\end{proof}
\begin{remark}
Let $\pi:(X,G)\rightarrow (Y,G)$ be a factor map between $G$-systems.
\begin{enumerate}[(1)]
  \item $D_e(X,G|\pi)=D(X,G|\pi)$ if one of the two values exists.
  \item If $(X,G)$ has a generating open cover, then there exists a relative entropy generating set (of the cover) $F$ such that $\overline{D}(F)=\overline{D}_e(X,G|\pi)=\overline{D}(X,G|\pi)$.
\end{enumerate}
\end{remark}
\section{Relative dimension tuples and dimension sets}\label{sec6}
In this section, we will localize relative entropy dimension to obtain the notion of relative dimension tuples and dimension sets. Let $\pi:(X,G)\rightarrow (Y,G)$ be a factor map between $G$-systems and $n\geq 2$. Denote by $\Delta_n(X)=\{(x_i)_1^n\in X^n:x_1=\cdots=x_n\}$ the $n$-th diagonal of $X$.
\subsection{Relative dimension tuples}\label{sec6.1}
\begin{defn}\label{d2}
Let $\pi:(X,G)\rightarrow (Y,G)$ be a factor map between $G$-systems and $(x_i)_1^n\in X^n\setminus\Delta_n(X)$. The relative entropy dimension of $(x_i)_1^n$ relevant to $\pi$ is
\[
\overline{D}(x_1,\cdots,x_n|\pi)=\lim\limits_{k\rightarrow \infty}\overline{D}(G,\mathcal{U}_k|\pi)\in[0,1],
\]
where $\mathcal{U}_k=\{X\setminus \overline{B(x_1,\frac{1}{k})},\cdots,X\setminus \overline{B(x_n,\frac{1}{k})}\}$. For $\alpha\in[0,1]$, when $\overline{D}(x_1,\cdots,x_n|\pi)=\alpha$, we will say $(x_i)_1^n$ is $\alpha$-relative dimension $n$-tuples relevant to $\pi$. The set of $\alpha$-relative dimension $n$-tuples relevant to $\pi$ of $(X,G)$ is denoted by $E_n^\alpha(X,G|\pi)$. Denote $G_n^\alpha(X,G|\pi)=\bigcup\limits_{\beta\geq\alpha}E_n^\beta(X,G|\pi)$.
\end{defn}
\begin{lem}\label{l2}
Let $\pi:(X,G)\rightarrow (Y,G)$ be a factor map between $G$-systems and $\mathcal{U}=\{U_1,\cdots,U_n\}\in\mathcal{C}_X^o$. Then there exists $x_i\in U_i^c$, $1\leq i\leq n$ such that $\overline{D}(x_1,\cdots,x_n|\pi)\geq\overline{D}(G,\mathcal{U}|\pi)$.
\end{lem}
\begin{proof}
We follow the idea in the proof of Proposition 2 in \cite{Bl}.

If $U_1^c$ is a singleton then we put $U_1^1=U_1$. If $U_1^c$ has at least two different points $y$ and $y'$, fix $\varepsilon_1>0$ with $\varepsilon_1\leq\frac{d(y_1,y_2)}{4}$, and construct a cover of $U_1^c$ by open balls with radius $\varepsilon_1$ centered in $U_1^c$; call it $\mathcal{A}$. Since $U_1^c$ is compact, there exist $A_1,A_2,\cdots,A_u\in\mathcal{U}$ such that $\bigcup_{i=1}^uA_i\supseteq U_1^c$. Write $F_i=U_1^c\cap\overline{A_i}$, $i=1,2,\cdots,u$. By the choice of $\varepsilon_1$, each closed set $F_i$ is a proper subset of $U_1^c$ with $diam(F_i)\leq \frac{1}{2}diam(U_1^c)$. Since $\{U_1,\cdots,U_n\}$ is coarser than $\bigvee_{i=1}^u\{F_i^c,U_2,\cdots,U_n\}$, we have
\[
\overline{D}(G,\mathcal{U}|\pi)\leq\overline{D}(G,\bigvee_{i=1}^u\{F_i^c,U_2,\cdots,U_n\}|\pi)=\max\limits_{1\leq i\leq u}\{\overline{D}(G,\{F_i^c,U_2,\cdots,U_n\}|\pi)\},
\]
where the last equality comes from Proposition \ref{p1}(3). Thus there exists $i_*\in\{1,\cdots,u\}$ such that $\overline{D}(G,\{F_{i_*}^c,U_2,\cdots,U_n\}|\pi)\geq\overline{D}(G,\mathcal{U}|\pi)$. We denote the set $F_{i_*}^c$ by $U_1^1$. Apply the same argument for $U_i$, and obtain $U_i^1$, $1<i\leq n$. Now we have found an open cover $\mathcal{U}_1=\{U_1^1,U_2^1,\cdots,U_n^1\}$, $1\leq i\leq n$ finer than $\mathcal{U}$ and having the property that \[\overline{D}(G,\mathcal{U}_1|\pi)\geq\overline{D}(G,\mathcal{U}|\pi)\mbox{ and }diam((U_i^1)^c)\leq\frac{1}{2}diam(U_i^c)\mbox{ for any }i=1,\cdots,n\]

Repeating the arguments, one can get $n$ decreasing sequences of non-empty closed sets $\{(U_1^j)^c\},\cdots,\{(U_n^j)^c\}$  and open cover $\mathcal{U}_j=\{U_1^j,U_2^j,\cdots,U_n^j\}$ such that
\begin{enumerate}[a)]
  \item $diam((U_i^{j+1})^c)\leq \frac{1}{2}diam((U_i^j)^c)$ for $i=1,\cdots,n$ and $j=0,1,2,\cdots$ (here let $U_i^0=U_i$, $i=1,\cdots,n$).
  \item $\overline{D}(G,\mathcal{U}_j|\pi)\geq \overline{D}(G,\mathcal{U}|\pi)$ for each $j=1,2,\cdots$.
\end{enumerate}
By the above a), we have $\bigcap_{j=1}^\infty(U_1^j)^c=\{x_1\}$, $\cdots$, $\bigcap_{j=1}^\infty(U_n^j)^c=\{x_n\}$ for some $x_1,\cdots,x_n$ in $X$. Moreover, since $\{U_1,\cdots,U_n\}$ covers $X$ and $x_i\in(U_i)^c$, we have $(x_i)_1^n\in X^n\setminus\Delta_n(X)$.

Finally, we show $\overline{D}(x_1,\cdots,x_n|\pi)\geq\overline{D}(G,\mathcal{U}|\pi)$. Since $(x_i)_1^n\in X^n\setminus\Delta_n(X)$, without loss of generality, we assume $x_1\neq x_n$. There exists $k_0\in\mathbb{N}$ such that $\overline{B(x_1,\frac{1}{k_0})}\cap\overline{B(x_n,\frac{1}{k_0})}=\emptyset$. For any $k\geq k_0$, there exists $j(k)\in\mathbb{N}$ large enough such that $(U_i^{j(k)})^c\subseteq B(x_i,\frac{1}{k})$ for $i=1,\cdots,n$, and thus
\[
\overline{D}(G,\{X\setminus \overline{B(x_1,\frac{1}{k})},\cdots,X\setminus \overline{B(x_n,\frac{1}{k})}\}|\pi)\geq\overline{D}(G,\mathcal{U}_{j(k)}|\pi)\geq \overline{D}(G,\mathcal{U}|\pi),
\]
where the last inequality comes from the above b). Hence we have
\[
\overline{D}(x_1,\cdots,x_n|\pi)=\lim\limits_{k\rightarrow \infty}\overline{D}(G,\{X\setminus \overline{B(x_1,\frac{1}{k})},\cdots,X\setminus \overline{B(x_n,\frac{1}{k})}\}|\pi)\geq \overline{D}(G,\mathcal{U}|\pi).
\]
\end{proof}
\begin{prop}\label{p6}
Let $\pi:(X,G)\rightarrow (Y,G)$ be a factor map between $G$-systems. Then
\begin{enumerate}[(1)]
  \item $G_n^\alpha(X,G|\pi)\cup \Delta_n(X)$ is a closed subset of $X^n$ for any $\alpha\geq0$.
  \item If $G$ is an abelian group, $G_n^\alpha(X,G|\pi)\cup \Delta_n(X)$ is $G$-invariant .
\end{enumerate}
\end{prop}
\begin{proof}
\begin{enumerate}[(1)]
\item Let $\{(x^i_1,\cdots,x^i_n)\}_{i=1}^{\infty}\subseteq G_n^{\alpha}(X,G|\pi)$ and $\lim\limits_{i\rightarrow +\infty}(x^i_1,\cdots,x^i_n)=(x_1,\cdots,x_n)$. If $(x_1,\cdots,x_n)\notin \Delta_nX$, then without loss of generality, we assume $x_1\neq x_n$. For any $k\in\mathbb{N}$ with $\frac{2}{k}<d(x_1,x_n)$, let \[\mathcal{U}_k=\{X\setminus \overline{B(x_1,\frac{1}{k}}),\cdots,X\setminus \overline{B(x_n,\frac{1}{k})}\}.\]
Then there exist $j\in\mathbb{N}$ and $k_j\in \mathbb{N}$ such that $\mathcal{U}_k\succeq\{X\setminus \overline{B(x_1^j,\frac{1}{k_j})},\cdots, X\setminus \overline{B(x^j_n,\frac{1}{k_j})}\}$. Therefore, we have
\[
\overline{D}(G,\mathcal{U}_k|\pi)\geq \overline{D}(G,\{X\setminus \overline{B(x_1^j,\frac{1}{k_j})},\cdots, X\setminus \overline{B(x^j_n,\frac{1}{k_j})}\}|\pi)\geq \overline{D}(x^j_1,\cdots,x^j_n|\pi)\geq \alpha.
\]
Let $k\rightarrow +\infty$, then we have $\overline{D}(x_1,\cdots,x_n|\pi)\geq \alpha$. Therefore $(x_1,\cdots,x_n)\in G^{\alpha}_n(G,X|\pi)$.
\item We now show $G^{\alpha}_n(G,X|\pi)\cup \Delta_X$ is $G$-invariant. Let $g\in G$ and let $(x_1,\cdots,x_n)\in G^{\alpha}_n(G,X|\pi)\cup \Delta_nX$. On the one hand, if $(x_1,\cdots,x_n)\in \Delta_nX$, then $(gx_1,\cdots,gx_n)\in \Delta_nX$. On the other hand, if $(x_1,\cdots,x_n)\in G^{\alpha}_n(G,X|\pi)$, then $\overline{D}(x_1,\cdots,x_n|\pi)\geq \alpha$. Then we have
\[
\begin{split}
\overline{D}(gx_1,\cdots,gx_n|\pi)&=\lim\limits_{k\rightarrow\infty}\overline{D}(G,\{X\setminus \overline{B(gx_1,\frac{1}{k})},\cdots,X\setminus \overline{B(gx_n,\frac{1}{k})}\}|\pi)\\
&=\lim\limits_{k\rightarrow\infty}\overline{D}(G,g^{-1}\{X\setminus \overline{B(gx_1,\frac{1}{k})},\cdots,X\setminus \overline{B(gx_n,\frac{1}{k})}\}|\pi)\\
&\geq\lim\limits_{m\rightarrow\infty}\overline{D}(G,\{X\setminus \overline{B(x_1,\frac{1}{m})},\cdots,X\setminus \overline{B(x_n,\frac{1}{m})}\}|\pi)=\overline{D}(x_1,\cdots,x_n|\pi)\geq \alpha,
\end{split}
\]
where the second equality comes from Proposition \ref{p1}(2).
Therefore, $(gx_1,\cdots,gx_n)\in G^{\alpha}_n(G,X|\pi)$, and then $G^{\alpha}_n(G,X|\pi)\cup \Delta_nX$ is $G$-invariant.
\end{enumerate}
\end{proof}
\begin{prop}\label{p7}
Let $(X,G)$, $(Y,G)$ and $(Z,G)$ be $G$-systems and $\pi$, $\pi_X$ and $\pi_Y$ are factor maps satisfying the following commutative diagram
\[
\xymatrix
{
(X,G)\small\ar[r]^\pi\small\ar[d]_{\pi_X}&(Y,G)\small\ar[dl]^{\pi_Y}\\
(Z,G),
}
\]

\begin{enumerate}[(1)]
  \item If $(x_1,\cdots,x_n)\in E_n^\alpha(X,G|\pi_X)$ and $y_i=\pi(x_i)$, $(y_i)_1^n\in Y^n\setminus \Delta_n(Y)$, then $\overline{D}(y_1,\cdots,y_n|\pi_Y)\geq\alpha$.
  \item If $(y_1,\cdots,y_n)\in E_n^\alpha(Y,G|\pi_Y)$, then there exists $(x_1,\cdots,x_n)\in E_n^\alpha(X,G|\pi_X)$ and $\pi(x_i)=y_i$, $i=1,\cdots,n$.
\end{enumerate}
\end{prop}
\begin{proof}
\begin{enumerate}[(1)]
  \item If $(x_1,\cdots,x_n)\in E_n^\alpha(X,G|\pi_X)$ and $y_i=\pi(x_i)$, $(y_i)_1^n\in Y^n\setminus \Delta_n(Y)$, then
  \[
  \begin{split}
  \overline{D}(y_1,\cdots,y_n|\pi_Y)=&\lim\limits_{k\rightarrow \infty}\overline{D}(G,\{Y\setminus \overline{B(y_1,\frac{1}{k})},\cdots,Y\setminus \overline{B(y_n,\frac{1}{k})}\}|\pi_Y)\\
  =&\lim\limits_{k\rightarrow \infty}\overline{D}(G,\pi^{-1}\{Y\setminus \overline{B(y_1,\frac{1}{k})},\cdots,Y\setminus \overline{B(y_n,\frac{1}{k})}\}|\pi_X)\\
  \geq&\lim\limits_{m\rightarrow \infty}\overline{D}(G,\{X\setminus \overline{B(x_1,\frac{1}{m})},\cdots,X\setminus \overline{B(x_n,\frac{1}{m})}\}|\pi_X)=\overline{D}(x_1,\cdots,x_n|\pi_X).
  \end{split}
  \]
  \item If $\alpha=0$, the conclusion holds obviously. Now assume $\alpha>0$ we let $(y_1,\cdots,y_n)\in E_n^\alpha(Y,G|\pi_Y)$. For any $k\in\mathbb{N}$ with $\frac{2}{k}<\max\limits_{1\leq i<j\leq n}d(y_i,y_j)$, let \[\mathcal{U}_k=\{Y\setminus \overline{B(y_1,\frac{1}{k})},\cdots,Y\setminus \overline{B(y_n,\frac{1}{k})}\}.\]
      Then by Lemma \ref{l2}, there exists
      \[
      (x_1^k,\cdots,x_n^k)\in \pi^{-1}(\overline{B(y_1,\frac{1}{k})})\times\cdots\times\pi^{-1}(\overline{B(y_n,\frac{1}{k})})
      \]
      with
      \[
      \overline{D}(x_1^k,\cdots,x_n^k|\pi_X)\geq\overline{D}(G,\pi^{-1}(\mathcal{U}_k)|\pi_Y\pi)=\overline{D}(G,\mathcal{U}_k|\pi_Y)\geq \alpha.
      \]
Hence $(x_1^k,\cdots,x_n^k)\in G^{\alpha}_n(X,G|\pi_X)$.
Take a subsequence $\{k_i\}$ such that
\[
\lim\limits_{i\rightarrow \infty}(x_1^{k_i},\cdots,x_n^{k_i})=(x_1,\cdots,x_n)
\]
for some $(x_1,\cdots,x_n)\in X^n$. Clearly $\pi(x_1)=y_1$, $\cdots$, $\pi(x_n)=y_n$, so $(x_1,\cdots,x_n)\notin\Delta_n(X)$. Now on the one hand, By Proposition \ref{p6}(1), we have $(x_1,\cdots,x_n)\in G^{\alpha}_n(X,G|\pi_X)$, that is, $\overline{D}(x_1,\cdots,x_n|\pi_X)\geq\alpha$. On the other hand, by (1), we have \[\overline{D}(x_1,\cdots,x_n|\pi_X)\leq\overline{D}(y_1,\cdots,y_n|\pi_Y)=\alpha.\] Therefore $\overline{D}(x_1,\cdots,x_n|\pi_X)=\alpha$,
i.e, $(x_1,\cdots,x_n)\in E_n^\alpha(X,G|\pi_X)$.
\end{enumerate}
\end{proof}
\subsection{Relative dimension sets and uniform relative dimension systems}\label{sec6.2}
\begin{defn}
Let $\pi:(X,G)\rightarrow (Y,G)$ be a factor map between $G$-systems. We call the subset $\{\alpha\geq0:E_n^\alpha(X,G|\pi)\neq\emptyset\}$ of $[0,1]$ the \textit{$n$-th relative dimension set} of $(X,G)$ relevant to $\pi$ and denote it by $\mathcal{D}_n(X,G|\pi)$. If $0\notin\mathcal{D}_n(X,G|\pi)$, we say $(X,G)$ is a system which has \textit{strictly positive $n$-th relative entropy dimension} relevant to extension $\pi$ . Let $\alpha\in(0,1]$, we call $(X,G)$ an \textit{$\alpha$-uniform $n$-th relative entropy dimension system} relevant to $\pi$ ($n$-th $\alpha$-u.r.d. system for short, if there is no confusion, we omit "n-th"), if $\mathcal{D}_n(X,G|\pi)=\{\alpha\}$ and call $(X,G)$ an \textit{$n$-th $\alpha^+$-relative dimension system} ($n$-th $\alpha^+$-r.d. system for short, if there is no confusion, we omit "n-th"), if $\mathcal{D}_n(X,G|\pi)\subseteq[\alpha,1]$. If $(X,G)$ is the trivial system, we let $\mathcal{D}_n(X,G|\pi)=\emptyset$.
\end{defn}
\begin{prop}
Let $\alpha\in(0,1]$. Then
\begin{enumerate}[(1)]
  \item A nontrivial $G$-system $(X,G)$ is a $n$-th $\alpha$-u.r.d. system relevant to $\pi$ if and only if $\overline{D}(G,\mathcal{U}|\pi)=\alpha$ for any open cover $\mathcal{U}$ of $X$ with $\mathcal{U}=\{U_1,\cdots,U_n\}$.
  \item In a commutative diagram
\[
\xymatrix
{
(X,G)\small\ar[r]^\pi\small\ar[d]_{\pi_X}&(Y,G)\small\ar[dl]^{\pi_Y}\\
(Z,G),
}
\]
if $(X,G)$ be an $\alpha$-u.r.d. system and $(Y,G)$ a nontrivial system, then $(Y,G)$ is also an $\alpha$-u.r.d. system.
  \item In the commutative diagram of (2), if $(X,G)$ be an $\alpha^+$-r.d. system relevant to $\pi_{X}$ and $(Y,G)$ a nontrivial system, then $(Y,G)$ is also an $\alpha^+$-r.d. system relevant to $\pi_{Y}$.
\end{enumerate}
\end{prop}
\begin{proof}
(2) and (3) come from Proposition \ref{p7}. We now show (1). Assume that $\overline{D}(G,\mathcal{U}|\pi)=\alpha$ for any open cover $\mathcal{U}$ of $X$ with $\mathcal{U}=\{U_1,\cdots,U_n\}$. Then by Definition \ref{d2}, for $(x_1,\cdots,x_n)\in X^n\setminus\Delta_n{X}$, we have $\overline{D}(x_1,\cdots,x_n|\pi)=\alpha$, that is, $(X,G)$ is $n$-th $\alpha$-u.r.d. system.

Conversely, suppose that $(X,G)$ is $n$-th $\alpha$-u.r.d. system. Let $\mathcal{U}=\{U_1,\cdots,U_n\}$ be an open cover of $X$. Then on the one hand, by Lemma \ref{l2}, there exists $x_i\in U_i^c$, $1\leq i\leq n$ such that $\overline{D}(G,\mathcal{U}|\pi)\leq\overline{D}(x_1,\cdots,x_n|\pi)=\alpha$. On the other hand, take $(x'_1,\cdots,x'_n)\in int(U_1^c)\times\cdots\times int(U_n^c)$, then $(x'_1,\cdots,x'_n)\notin \Delta_n(X)$ and $\overline{D}(G,\mathcal{U}|\pi)\geq\overline{D}(x'_1,\cdots,x'_n|\pi)=\alpha$. Hence $\overline{D}(G,\mathcal{U}|\pi)=\alpha$.
\end{proof}
\section{A relative disjointness theorem via entropy dimensions}\label{sec7}

Let $\pi_X:(X,G)\rightarrow (Z,G)$ and $\pi_Y:(Y,G)\rightarrow (Z,G)$ be two factor maps, and $\pi_1:X\times Y\rightarrow X$, $\pi_2:X\times Y\rightarrow Y$ be two projections. $J\subseteq X\times Y$ is called a \textit{joining} of $(X,G)$ and $(Y,G)$ over $(Z,G)$, if $J$ is closed and $G$-invariant with
\[
\pi_1(J)=X,\pi_2(J)=Y \mbox{ and }\pi_X\times\pi_Y(J)=\Delta_2(Z)(\mbox{i.e. }\pi_X\pi_1=\pi_Y\pi_2).
\]
Define \[X\times_Z Y=\bigcup\limits_{z\in Z}\pi_X^{-1}(z)\times\pi_Y^{-1}(z).\]

Clearly, $X\times_ZY$ is a joining of $(X,G)$ and $(Y,G)$ over $(Z,G)$. A joining $J$ of $(X,G)$ and $(Y,G)$ over $(Z,G)$ is said to be \textit{proper} if $J\neq X\times_ZY$. We say that $(X,G)$ and $(Y,G)$ are \textit{disjoint} over $(Z,G)$, if $X\times_ZY$ contains no proper joining of $(X,G)$ and $(Y,G)$ over $(Z,G)$.

Let $\pi:(X,G)\rightarrow (Y,G)$ be a factor map between two $G$-systems. We say $\pi$ is \textit{minimal} if $X$ is the only closed and $G$-invariant subset with $\pi$-image $Y$, and that $\pi$ is said to be \textit{open} if the $\pi$-image of any open subset is open.

The following result is a relative version of \cite[Theorem 4.12]{DHP} which generalizes the theorem that uniformly positive entropy systems are disjoint from minimal entropy zero systems \cite{Bl}.
\begin{thm}\label{t2}
Let $\pi_X:(X,G)\rightarrow (Z,G)$ and $\pi_Y:(Y,G)\rightarrow (Z,G)$ be two factor maps with $\pi_X$ open and $\pi_Y$ minimal. If for any $n\geq 2$, $\mathcal{D}_n(G,X|\pi_X)>\mathcal{D}_n(G,Y|\pi_Y)$ (i.e. for any $\alpha\in\mathcal{D}_n(G,X|\pi_X)$ and $\beta\in\mathcal{D}_n(G,Y|\pi_Y)$, $\alpha>\beta$), then $(X,G)$ and $(Y,G)$ are disjoint over $(Z,G)$.
\end{thm}
\begin{proof}
Let $J\subseteq X\times_ZY$ be any given joining of $(X,G)$ and $(Y,G)$ over $(Z,G)$.

Let $J(x)=\{y\in Y:(x,y)\in J\}$ for each $x\in X$. For each $z\in Z$, we denote \[J^*(z)=\bigcap\limits_{x\in\pi_X^{-1}(z)}J(x).\]
First we claim that $J^*(z)\neq\emptyset$ for $z\in Z$. Assume the contrary that there exists $\{x_i\}_1^n\subseteq \pi_X^{-1}(z)$ with $\bigcap\limits_{i=1}^nJ(x_i)=\emptyset$. Clearly, $(x_1,\cdots,x_n)\notin\Delta_n(X)$. Let $\alpha=\overline{D}(x_1,\cdots,x_n|\pi_X)$. Then $\alpha\in\mathcal{D}_n(G,X|\pi_X)$. By Proposition \ref{p7}(2) for $\pi_1:J\rightarrow X$, there exists $y_i\in Y$ such that
\[
(x_i,y_i)\in J \mbox{ and }\overline{D}((x_1,y_1),\cdots,(x_n,y_n)|\pi_X\pi_1)=\alpha.
\]
Since $\bigcap\limits_{i=1}^nJ(x_i)=\emptyset$, we have $(y_1,\cdots,y_n)\notin\Delta_n(Y)$. Then by Proposition \ref{p7}(1) for $\pi_2:J\rightarrow Y$ and $\pi_X\pi_1=\pi_Y\pi_2$, we have
\[
\overline{D}(y_1,\cdots,y_n|\pi_Y)\geq\alpha,
\]
which is a contradiction to the assumption that $\mathcal{D}_n(G,X|\pi_X)>\mathcal{D}_n(G,Y|\pi_Y)$.

In order to prove that $(X,G)$ and $(Y,G)$ are disjoint over $(Z,G)$, it is sufficient to prove $J=X\times_ZY$. Put
\[
J'=\bigcup\limits_{z\in Z}\pi_X^{-1}(z)\times J^*(z)\subseteq J\subseteq X\times_ZY.
\]

Obviously, $J'$ is $G$-invariant. Now we show $J'$ is closed. Let $z_n\in Z$ and $y_n\in J^*(z_n)$ such that $y_n\rightarrow y_0$ and $z_n\rightarrow z_0$. Since $\pi_X$ is open, the map $z\mapsto\pi_X^{-1}(z)$ is continuous. For any $x_0\in \pi_X^{-1}(z_0)$, the continuity implies that there exists $x_n\in\pi^{-1}_X(z_n)$ with $x_n\rightarrow x_0$. Since $y_n\in J(x_n)$ and the map $x\mapsto J(x)$ is upper semi-continuous, $y_0=\lim\limits_ny_n\in J(x_0)$. Since $x_0$ is arbitrary, $y_0\in J^*(z_0)$. This proves $J'$ is closed. Then by the previous claim, $J^*(z)\neq\emptyset$ for $z\in Z$. Moreover, $\pi_1(J')=\bigcup\limits_{z\in Z}\pi_X^{-1}(z)=X$.

Note that $\pi_2(J')$ is closed and $G$-invariant subset of $Y$ and $\pi_Y\pi_2(J')=\pi_X\pi_1(J')=\pi_X(X)=Z$. We deduce that $\pi_2(J')=Y$ by the minimality of $\pi_Y$. Since $\pi_2(J')=\bigcup\limits_{z\in Z}J^*(z)$ and $J^*(z)\subseteq\pi_Y^{-1}(z)$, we get $J^*(z)=\pi_Y^{-1}(z)$ for all $z\in Z$. This implies $J'=J=X\times_ZY$, i.e. $(X,G)$ and $(Y,G)$ are disjoint over $(Z,G)$.
\end{proof}
\section*{\textbf{Acknowledgments}}%

The first author was partly supported by the Educational Research Project
for Young and Middle-aged Teachers of Fujian Province (Grant No. JAT200045).
\section*{\textbf{References}}%

\end{document}